\documentclass{amsart}

\newtheorem{theorem}{Theorem}[section]

\newtheorem{propo}{Proposition}[section]
\newtheorem{corollary}[theorem]{Corollary}
\theoremstyle{definition}
\newtheorem{definition}[theorem]{Definition}

\theoremstyle{remark}
\newtheorem{remark}[theorem]{Remark}
\numberwithin{equation}{section}

\begin{document}
	
	\title{ Rigidity results on $\rho$--Einstein solitons with zero scalar
		 curvature }
	
	\author{Romildo Pina}
	\curraddr{Instituto de Matem\'atica e Estat\'istica, Universidade Federal de Goi\'as, Goi\^ania, Brasil, 74001-970}
	\email{romildo@ufg.br}
	
	\author{Ilton Menezes}
	\email{iltomenezesufg@gmail.com}
		\author{Lucyjane Silva}
	\email{lucyjanedealmeida@gmail.com}


	
	\begin{abstract}
		In this paper we show that a $\rho$-Einstein solitons conformal to a 
		pseudo-Euclidean space, invariant under the action of the pseudo-orthogonal group with zero scalar curvature is stady and consequently flat. How
		application of the results obtained we present an explicit example for a the question proposed by Kazdan in \cite{Ka}. 
		
	\end{abstract}
	
	\maketitle
	
\section{Introduction and main statements} 

In this paper, we study two related problems. The first problem is on the existence of $\rho$--Einstein sotiton with scalar curvature $K_{\bar{g}}=0$. Besides that we present some rigity results.

The second problem consists in find all metrics that are conformal to the pseudo Euclidean metrics, with zero scalar curvature, which are invariant under the action of the pseudo-orthogonal group. This provides explicit solutions to Yamabe's problem in the non-compact case. In the Riemannian case under some additional assumptions, all metrics obtained are complete. As application of this results we obtain a family of complete metrics in $\mathbb{R}^{n}\setminus\{0\}$ with scalar curvature positive, negative and zero, presenting an explicit example for a question proposed by
Kazdan in \cite{Ka}.

In 1982, R. Hamilton introduced a nonlinear evolution equation for Riemannian metrics with the aim of finding canonical metrics on manifolds (see \cite{SB} or \cite{Ha}). This evolution equation is known as the Ricci flow, and it has since been used widely and with great success, most notably in Perelman's solution of the Poincaré conjecture. Furthermore, several convergence theorems have been established. One important aspect in the treatment of the Ricci flow is the study of Ricci solitons, which generate self-similar solutions to the flow and often arise as singularity models.

Given a semi-Riemannian manifold $(M^{n},g)$, $n\geq 3$, we say that $(M,g)$ is a gradient Ricci soliton if there exists a differentiable function $h:M\longrightarrow \mathbb{R}$ (called the potential function) such that
\begin{equation}\label{000}
\mbox{Ric}_g+\mbox{Hess}_g(h)=\lambda g, \qquad \lambda\in\mathbb{R},
\end{equation}
where Ric$_g$ is the Ricci tensor, Hess$_g(h)$ is the Hessian of $h$ 
with respect to the metric $g$, and $\lambda$ is a real number. We say that a gradient Ricci soliton is {\em shrinking, steady, or expanding} if 
$\lambda>0$, $\lambda=0$, or $\lambda<0$, respectively.
Bryant \cite{BRYANT} proved that there exists a 
complete, steady,  gradient Ricci soliton that is spherically symmetric for any 
$n\geq 3$, which is known as Bryant's soliton. In the bi-dimensional case an analogous nontrivial rotationally symmetric solution was obtained explicitly, and is known as the Hamilton cigar.
Recently  Cao-Chen \cite{CAOCHEN} showed that any complete, steady, gradient Ricci 
soliton, locally conformally flat, up to homothety, is either flat or isometric to 
the Bryant's soliton.  The results obtained in  \cite{CAOCHEN}  were extended to  bach -- flat gradient  steady Ricci solitons (see \cite{CCM}). Complete, conformally flat shrinking gradient 
solitons have been characterized as being quotients of
$\mathbb{R}^n$, $\mathbb{S}^n$ or $R\times \mathbb{S}^{n-1}$ (see \cite{FG}). In the case of steady gradient Ricci solitons, \cite{RO}  provides all such solutions when the metric is  conformal to an n-dimensional pseudo-Euclidean space and invariant under the action of an $(n-1)-$dimensional translation group. 

Motivated by the notion of Ricci solitons on a semi-Riemannian manifold $(M^{n},g)$, $n\geq3$, it is natural to consider geometric flows of the following type: 
\begin{equation}\label{flow}
\frac{\partial}{\partial t}g(t)=-2(Ric-\rho Rg)
\end{equation}
for $\rho\in\mathbb{R}$, $\rho\neq0$, as in \cite{BO} and \cite{CA}. We call these the  Ricci-Bourguignon flows. We notice that  short time existence for the geometric flows described in \eqref{flow} is provided
in (\cite{CCD}).  Associated to the flows, we have the following notion of gradient $\rho$-Einstein solitons, which generate self-similar solutions:

\begin{definition} Let $\left(M^{n},g\right), n\geq3$, be a Riemannian manifold and let $\rho\in \mathbb{R},\rho\neq0$. We say that $(M^{n},g)$ is a gradient $\rho-$Einstein soliton if there exists a smooth function $h:M\longrightarrow\mathbb{R}$, such that the metric $g$ satisfies the equations 
	\begin{align}\label{def. 1}
	Ric_{g}+Hess_{g}h=\rho K_{g}g +\lambda g
	\end{align}
	for some constant $\lambda\in\mathbb{R}$, where $K_{g}$ is the scalar curvature of the metric $g$.
\end{definition}

A $\rho$-Einstein soliton is said to be shrinking, steady, or expanding  if $\lambda>0$, $\lambda=0$, or $\lambda<0$, respectively. Furthermore, a $\rho$-Einstein solitons is said to be a gradient Einstein soliton, gradient traceless Ricci soliton, and gradient Schouten soliton if $\rho=\frac{1}{2}$, $\rho=\frac{1}{n}$, and $\rho=\frac{1}{2(n-1)}$, respectively.

The gradient $\rho-$Einstein solitons equation \eqref{def. 1} links geometric information about the curvature
of the manifold through the Ricci tensor and the geometry of the level sets of the
potential function by means of their second fundamental form. Hence, classifying
gradient $\rho-$Einstein solitons under some curvature conditions is a natural problem. The $\rho-$Einstein solitons were investigated by Catino and Mazzieri
in \cite{CA}, they obtained important rigidity results, proving that every compact gradient Einstein, Schouten, or traceless Ricci soliton is trivial. In addition, they proved that every complete gradient steady Schouten soliton is trivial, hence Ricci flat.

Gradient Ricci solitons with constant scalar curvature were investigated by Petersen
and Wylie in \cite{P}, they proved that: If a non-steady gradient Ricci soliton has constant scalar curvature $K_{g}$, then it is
bounded as $0 \leq K_{g} \leq n\lambda$ in the shrinking case, and $n\lambda \leq K_{g} \leq 0$ in the expanding case. Fern\'andez-L\'opez and Garcia-R\'io in \cite{GR} improved this result proving that: If an n--dimensional complete gradient Ricci soliton with
constant scalar curvature $K_{g}$, then $K_{g}$ must be a multiple of $\lambda$.

In \cite{MP} the authors considered a $\rho$--Einstein solitons that are conformal to a pseudo Euclidean space and invariant under the action of the pseudo-orthogonal group.
They provide all the solutions for the gradient Schouten soliton case. Moreover,
they proved that if a gradient Schouten soliton is both complete, conformal to a
Euclidean metric, and rotationally symmetric, then it is isometric to $\mathbb{R}\times\mathbb{S}^{n-1}$.

In \cite{Ma} the autors used the variational method to study the existence problem of metrics with constant scalar curvature on complete non-compact Riemannian manifolds. The assumptions of results are motivated from question  in the work of Kazdan \cite{Ka}. The question
is that if M has complete metrics $g_{+}$ and $g_{-}$ with positive (respectively negative) scalar curvature, is there one with zero scalar curvature? With several additional hypotheses o autor provide an answer to the question posed by Kazdan, more details see \cite{MA}. 

We studied the equation \eqref{def. 1} in semi-Riemannian manifolds with scalar curvature constante. We consider gradient $\rho$-Einstein solitons conformal to a pseudo-Euclidean space, which are invariant under the action of the pseudo-orthogonal group. More precisely, let $(\mathbb{R}^{n},g)$ be the standard pseudo-Euclidean space with metric $g$ and coordinates $(x_{1},...,x_{n})$, with $g_{ij}=\delta_{ij}\varepsilon_{i}$, $1\leq i,j \leq n$, where $\delta_{ij}$ is the Kronecker delta, and $\varepsilon_{i}=\pm1$. Let $r=\sum_{i=1}^{n}\varepsilon_{i}x_{i}^{2}$ be a basic invariant for an $(n-1)-$dimensional
pseudo-orthogonal group. The main goal of this paper is to present in the Riemannian case a family of complete metrics and some results of rigidity
 on a large class of noncompact semi-Riemannian manifolds in the case where the scalar curvature is zero. In the Riemannian case the same results hold. 

We initially find a system of differential equations, such that the functions $h$ and $\psi$ must satisfy, so that the metric $\bar{g}= g/\psi^{2}$ satisfies \eqref{def. 1} (see Theorem \ref{theorem 1}). Note that if the solutions are invariant under the action of the pseudo-orthogonal group, the system of partial differential equations given in Theorem \ref{theorem 1} can be transformed into a system of ordinary differential equations (see Corollary \ref{theorem 2}). In the Theorem \ref{coro1} we found all metrics that are conformal to the pseudo Euclidean metrics, with zero scalar curvature, which are invariant under the action of the pseudo-orthogonal group. As a consequence of the Theorem \ref{coro1}, we obtain the Corollary \ref{coro2} we constuct a family of complete metrics with zero scalar curvature. We present  results of rigidity on  gradient $\rho$- Einstein soliton with scalar curvature zero ( Theorem\ref{theorem 5}, Corollary \ref{theorem 6} and Corollary \ref{theorem 7}).

 In the Proposition \ref{prop} we construct a family of complete metrics with zero scalar curvature. In the Corollary \ref{coro3} we constuct an explicit example for Kazdan's question.

In what follows, we state our main results. We denote the second order derivative of $\psi$ and $h$ by $\psi_{,x_{i}x_{j}}$ and $h_{,x_{i}x_{j}}$, respectively, with respect to $x_{i}x_{j}$.
\begin{theorem}\label{theorem 1}
	Let $\left(\mathbb{R}^{n},g\right)$,$n\geq 3$, be a pseudo-Euclidean space with  coordinates $x=\left(x_{1},...,x_{n}\right)$ and metric components $g_{ij}=\delta_{ij}\varepsilon_{i}$,  $1\leq i,j \leq n$, where $\varepsilon_{i}=\pm 1$. Consider a smooth function $h:\mathbb{R}^{n}\longrightarrow \mathbb{R}$. There exists a metric $\bar{g}=\frac{1}{\psi^{2}}g$ such that $\left(\mathbb{R}^{n},\bar{g}\right)$ is a gradient $\rho$-Einstein soliton with $h$ as a potential function if, and only if, the functions $\psi$, and $h$ satisfy
	\begin{align}\label{01}
	(n-2)\psi_{,x_{i}x_{j}}+\psi h_{,x_{i}x_{j}}+\psi_{,x_{i}}h_{,x_{j}}+\psi_{,x_{j}}h_{,x_{i}}=0,\hspace{0,2cm} i\neq j,
	\end{align}
and
	\begin{center}
		\begin{align}\label{2}
		\psi\left[\left(n-2\right)\psi_{,x_{i}x_{i}}+\psi h_{,x_{i}x_{i}}+2\psi_{,x_{i}}h_{,x_{i}} \right]
		\end{align}
		\begin{align*}
		+\varepsilon_{i}\sum_{k=1}^{n}\varepsilon_{k}\left[\left(n-1\right)\left(\rho n \psi_{,x_{k}}^{2}-2\rho \psi \psi_{,x_{k}x_{k}}-\psi_{,x_{k}}^{2}\right)-\psi \psi_{,x_{k}}h_{,x_{k}}+\psi \psi_{,x_{k}x_{k}}\right]=\lambda \varepsilon_{i}, \hspace{0,2cm} i=j.
		\end{align*}
	\end{center}
	
\end{theorem}

Our objective is to determine solutions of the system \eqref{01}, \eqref{2} of the form $\psi(r)$ and $h(r)$, where $r=\sum_{i=1}^{n}\varepsilon_{i}x_{i}^{2}$. The following theorem reduces the system of partial differential equations \eqref{01} and \eqref{2} into an system of ordinary differential equations.
\vspace{.2in}

\begin{corollary}\label{theorem 2}
	Let $\left(\mathbb{R}^{n},g\right)$,$n\geq 3$, be a pseudo-Euclidean space with  coordinates $x=\left(x_{1},...,x_{n}\right)$ and metric components $g_{ij}=\delta_{ij}\varepsilon_{i}$, $1\leq i,j \leq n$, where $\varepsilon_{i}=\pm 1$. Consider smooth functions $\psi(r)$ and $h(r)$ with $r=\sum_{k=1}^{n}\varepsilon_{k}x_{k}^{2}$. Then there exists a metric $\bar{g}=\frac{1}{\psi^{2}}g$ such that $\left(\mathbb{R}^{n},\bar{g}\right)$ is a gradient $\rho$-Einstein soliton with $h$ as a potential function if, and only if, the functions $\psi$ and $h$ satisfy
	\vspace{12pt}
	\begin{equation}\label{3}
	(n-2)\psi''+\psi h''+2\psi'h'=0,
	\end{equation}
	and 
	\begin{center}
		\begin{equation}\label{4}
		2\psi\left[\left(n-2\right)\psi'+\psi h'\right]+2n[1-2(n-1)\rho]\psi \psi'
		\end{equation}
		\begin{equation*}
		+4r\left\{\left(n-1\right)\left[\left(\rho n-1\right)(\psi')^{2}-2\rho\psi \psi''\right]-\psi \psi'h'+\psi \psi''\right\}=\lambda.
		\end{equation*}
	\end{center}
\end{corollary}
  
The next we found all metrics that are conformal to the pseudo Euclidean metrics, with zero scalar curvature, which are invariant under the action of the pseudo-orthogonal group.

\begin{theorem}\label{coro1}
	Let $( \mathbb{R}^n, g)$ be a
	pseudo-Euclidean space, $n\geq 3$, with coordinates
	$x=(x_1,\cdots, x_n)$ and $g_{ij}=\delta_{ij}\varepsilon_i$, $1\leq i,j\leq n$, where $\varepsilon_i=\pm1$.
	Consider $\bar{g}=\frac{1}{\psi(r)^{2}}g$ where 	$r=\sum\limits_{k=1}^n\varepsilon_kx_{k}^2$. Then $\bar{g}$ have scalar curvature $K_{\bar{g}}=0$, if and only if

	\begin{equation}
	\psi(r)= \frac{k_2r}{\left(1+Ar^{\frac{n-2}{2}}\right)^{\frac{2}{n-2}}},
	\end{equation}
	where  $A, k_2 \in \mathbb{R}$ with $k_2>0$. If $A \geq 0$ the metric $\bar{g}$ is defined in $ \mathbb{R}^n\setminus\{0\}$. If $A<0$ the set of singularity points of $\bar{g}$ consist of the origin and a sphere $(n-1)$--dimensional, with center at the origin and radius $R=\sqrt{(\frac{-1}{A})^{\frac{2}{n-2}}}$.   
\end{theorem}
\begin{remark}
If $(\mathbb{R}^n, g)$ is the Euclidean space, then we find in the Theorem \ref{coro1} all metrics conformal  to $g$ and spherically symmetrical with zero scalar curvature. This provides explicit solutions to Yamabe's problem in the non-compact case.
\end{remark}

In \cite{CO}, the authors showed that $\{\mathbb{R}^n\setminus\{0\}, \bar{g}=\frac{1}{\varphi^2}g_{0}, \varphi(r)=\sqrt{r}\}$ 
is a complete Riemannian manifold and isometric at $\mathbb{S}^{n-1}\times\mathbb{R}$. As a consequence of Theorem \ref{coro1} together with this fact, we obtain the following result:
\begin{corollary}\label{coro2}
	Let $( \mathbb{R}^n, g)$ be a
	Euclidean space, $n\geq 3$, with coordinates
	$x=(x_1,\cdots, x_n)$ and $\left(g_{0}\right)_{ij}=\delta_{ij}$, $1\leq i,j\leq n$.
	Consider $\bar{g}=\frac{1}{\psi(r)^{2}}g_{0}$ where 	$r=\sum\limits_{k=1}^nx_{k}^2$. The metrics obtained in the Theorem \ref{coro1} are complete whenever $A>0$.
\end{corollary}

As an consequence of the Theorem \ref{coro1}, we get the following rigidity results.
\begin{theorem}\label{theorem 5}
	Let $( \mathbb{R}^n, g)$ be a
	pseudo-Euclidean space, $n\geq 3$, with coordinates
	$x=(x_1,\cdots, x_n)$ and $g_{ij}=\delta_{ij}\varepsilon_i$, $1\leq i,j\leq n$, where $\varepsilon_i=\pm1$.
	Consider $\left(\mathbb{R}^{n},\bar{g}\right)$, $\bar{g}=\frac{1}{\psi^{2}}g$  a  $\rho$--Einstein sotiton with scalar curvature $K_{\bar{g}}=0$, where  $\psi(r)$ and $h(r)$  smooth functions, 
	$r=\sum\limits_{k=1}^n\varepsilon_kx_{k}^2$ and $h$ as a potential function.	Then $\lambda=0$, that is $\left(\mathbb{R}^{n},\bar{g}\right)$ is steady.
\end{theorem}

\begin{corollary}\label{theorem 6}
	Let $( \mathbb{R}^n, g)$ be a
	pseudo-Euclidean space, $n\geq 3$, with coordinates
	$x=(x_1,\cdots, x_n)$ and $g_{ij}=\delta_{ij}\varepsilon_i$, $1\leq i,j\leq n$, where $\varepsilon_i=\pm1$. Then  $\left(\mathbb{R}^{n},\bar{g}\right)$, $\bar{g}=\frac{1}{\psi^{2}}g$ is  a  $\rho$--Einstein sotiton steady with scalar curvature $K_{\bar{g}}=0$, where  $\psi(r)$ and $h(r)$  smooth functions, 
	$r=\sum\limits_{k=1}^n\varepsilon_kx_{k}^2$ and $h$ as a potential function, if and only if, $\left(\mathbb{R}^{n},\bar{g}\right)$ is flat. 
			
\end{corollary}

As a consequence of the previous results, we have the following result in the Riemannian case.
\begin{corollary}\label{theorem 7}
Let $( \mathbb{M}^n, \bar{g})$ be $n\geq 3$ a  $\rho$--Einstein sotiton, Riemannian,locally conformally  flat and rotationally symmetric with zero scalar curvature.Then 
$(\mathbb{M}^n, \bar{g})$ is necessarily steady. Besides that $(\mathbb{M}^n, \bar{g})$ is flat. 
\end{corollary}

\begin{remark}	
	These results hold for $\rho=0$ and therefore they are extended to the Ricci solitons gradients,  proving that a Ricci soliton gradient, conformal  to the Euclidean space and spherically symmetrical with zero scalar curvature is necessarily steady and consequently flat. 
\end{remark}

\begin{remark}	
 As a consequently of the results obtained let's make an aplication giving a positive answer for a question proposed by Kazdan in \cite{Ka}, as follows:
 
 If M has complete metrics $g_+$ and $g_ -$ with positive (respectively, negative)
scalar curvature, is there one with zero scalar curvature?  kazdan mostrou em  \cite{KW},  que no caso compacto the answer is
"yes".

We built in $(\mathbb{R}^n\setminus\{0\})$ complete metrics with positive, negative and zero scalar curvature, respectively.
 
\end{remark}

\begin{propo}\label{prop}
	Let $( \mathbb{R}^n, g)$ be a
	Euclidean space, $n\geq 3$, with coordinates
	$x=(x_1,\cdots, x_n)$ and $\left(g_{0}\right)_{ij}=\delta_{ij}$, $1\leq i,j\leq n$. Consider $g=\frac{1}{\varphi(r)^{2}}g_{0}$ where 	$r=\sum\limits_{k=1}^nx_{k}^2$. If $\varphi(r)=re^{-\left(1+r^{\frac{2}{n-2}}\right)^{\frac{n-2}{2}}}$, then the metric $g$ on $\mathbb{R}^n$ is complete with scalar curvature negative given by 
	\begin{equation*}
		K_{g}=h(r)\left[(n-2)\left(1+r^{\frac{n-2}{2}}\right)^{\frac{2}{n-2}}+2(n-1)r^{\frac{2-n}{2}}+(n+2)\right],
	\end{equation*}
	where $h(r)=-\frac{4(n-1)r^{n-1}\left(1+r^{\frac{n-2}{2}}\right)^{\frac{2(3-n)}{n-2}}}{e^{2\left(1+r^{\frac{n-2}{2}}\right)^{\frac{2}{n-2}}}}$.
\end{propo}

In the next result we construct an explicit example in Riemannian manifolds for the question left by Kazdan \cite{Ka}.
\begin{corollary}\label{coro3} 
	Note that $\{\mathbb{R}^n\setminus\{0\}, \bar{g}=\frac{1}{\varphi^2}g_{0}, \varphi(r)=\sqrt{r}\}$ 
	is a complete Riemannian manifold with scalar curvature positive and $\left\{\mathbb{R}^n\setminus\{0\}, \bar{g}=\frac{1}{\varphi_1^2}g_{0}, \varphi_1(r)=re^{-\left(1+r^{\frac{2}{n-2}}\right)^{\frac{n-2}{2}}}\right\}$ 
	is a complete Riemannian manifold with scalar curvature negative, exists a complete metric of scalar curvature zero.
\end{corollary}

\section{Proofs of the main results}\label{sec2}

\begin{proof}
	Proof of Theorem \ref{theorem 1}. It is well known (see, e.g., \cite{RO}) that if $\bar{g}=\frac{g}{\psi^{2}}$, then 
	
	\begin{align*}
	Ric_{\bar{g}}=\frac{1}{\psi^{2}} \{(n-2)\psi Hess_{g}(\psi)+[\psi\Delta_{g}\psi-(n-1)|\nabla_{g}\psi^{2}]g \}
	\end{align*}
	and 
	\begin{align*}
	\bar{R}=(n-1)\left(2\psi\Delta_{g}\psi-n|\nabla_{g}\psi|^{2}\right).
	\end{align*}
	
	Hence, the equation
	\begin{align*}
	Ric_{\bar{g}}+Hess_{\bar{g}}(h)=\rho \bar{R}\bar{g}+\lambda\bar{g},
	\end{align*}
	is equivalent to 
	\begin{align}\label{A}
	\frac{1}{\psi^{2}} \{(n-2)\psi Hess_{g}(\psi)_{ij}+[\psi\Delta_{g}\psi-(n-1)|\nabla_{g}\psi|^{2}]\delta_{ij}\varepsilon_{i}\}+Hess_{\bar{g}}(h)_{ij}
	\end{align}
	\begin{equation*}
	=\left[\rho(n-1)(2\psi\Delta_{g}\psi-n|\nabla_{g}\psi|^{2})+\lambda\right]\frac{1}{\psi^{2}}\delta_{ij}\varepsilon_{i}.
	\end{equation*}
	
	Recall that, 
	\begin{equation*}
	Hess_{\bar{g}}(h)_{ij}=h_{,x_{i}x_{j}}-\sum_{k=1}^{n}\bar{\Gamma}_{ij}^{k}h_{,x_{k}}
	\end{equation*}
	where $\bar{\Gamma}_{ij}^{k}$ are the Christoffel symbols of the metric $\bar{g}$. For a distinct $i,j,k$, we have 
	\begin{equation*}
	\bar{\Gamma}_{ij}^{k}=0,\hspace{0.5cm}\bar{\Gamma}_{ij}^{i}=-\frac{\psi_{,x_{j}}}{\psi},\hspace{0.5cm}\bar{\Gamma}_{ii}^{k}=\varepsilon_{i}\varepsilon_{k}\frac{\psi_{,x_{k}}}{\psi},\hspace{0.5cm}\bar{\Gamma}_{ii}^{i}=-\frac{\psi_{,x_{i}}}{\psi},
	\end{equation*}
	therefore,
	
	\begin{equation}\label{10}
	Hess_{\bar{g}}(h)_{ij}=
	h_{,x_{i}x_{j}}+\frac{\psi_{,x_{j}}h_{,x_{i}}}{\psi}+\frac{\psi_{,x_{i}}h_{,x_{j}}}{\psi}, \hspace{0,2cm} i\neq j.
	\end{equation}
	
	Similarly, by considering $i=j$, we have 
	
	\begin{equation}\label{11}
	Hess_{\bar{g}}(h)_{ii}=
	h_{,,x_{i}x_{i}}+\frac{2\psi_{,x_{i}}h_{,x_{i}}}{\psi}-\varepsilon_{i}\sum_{k=1}^{n}\varepsilon_{k}\frac{\psi_{,x_{k}}h_{,x_{k}}}{\psi}.
	\end{equation}
	
	However, we note that
	\begin{align}\label{14}
	|\nabla_{g}\psi|^{2}=\sum_{k=1}^{n}\varepsilon_{k}\left(\frac{\partial\psi}{\partial x_{k}}\right)^{2}, \hspace{0.5cm} \Delta_{g}\psi=\sum_{k=1}^{n}\varepsilon_{k}\psi_{,x_{k}x_{k}},\hspace{0.5cm} Hess_{g}(\psi)_{ij}=\psi_{,x_{i}x_{j}}.
	\end{align}
	
	If $i\neq j$ in \eqref{A}, we obtain 
	\begin{align}\label{15}
	(n-2) \frac{Hess_{g}(\psi)_{ij}}{\psi}+Hess_{\bar{g}}(h)_{ij}=0.
	\end{align}
	
	Substituting the expressions found in \eqref{10}, and \eqref{14} into \eqref{15}, we obtain 
	
	\begin{align*}
	(n-2)\psi_{,x_{i}x_{j}}+\psi h_{,x_{i}x_{j}}+\psi_{,x_{i}}h_{,x_{j}}+\psi_{,x_{j}}h_{,x_{i}}=0,\hspace{0,2cm} i\neq j. 
	\end{align*}
	
	Similarly, if $i=j$ in \eqref{A}, we have 
	\begin{align}\label{13}
	(n-2)\psi Hess_{g}(\psi)_{ii}+\psi\Delta_{g}\psi\varepsilon_{i}-(n-1)|\nabla_{g}\psi|^{2}\varepsilon_{i}+\psi^{2}Hess_{\bar{g}}(h)_{ii}
	\end{align}
	\begin{align*}
	=2(n-1)\rho\Delta_{g}\psi \varepsilon_{i}-n(n-1)\rho|\nabla_{g}\psi|^{2}\varepsilon_{i}+\lambda\varepsilon_{i}.
	\end{align*}
	
	Substituting the expressions found in \eqref{11}, and \eqref{14} into \eqref{13}, we obtain 
	
	\begin{align*}
	\psi\left[\left(n-2\right)\psi_{,x_{i}x_{i}}+\psi h_{,x_{i}x_{i}}+2\psi_{,x_{i}}h_{,x_{i}} \right]
	\end{align*}
	\begin{align*}
	+\varepsilon_{i}\sum_{k=1}^{n}\varepsilon_{k}\left[\left(n-1\right)\left(\rho n \psi_{,x_{k}}^{2}-2\rho \psi \psi_{,x_{k}x_{k}}-\psi_{,x_{k}}^{2}\right)-\psi \psi_{,x_{k}}h_{,x_{k}}+\psi \psi_{,x_{k}x_{k}}\right]=\lambda \varepsilon_{i}.
	\end{align*}
	
	This concludes the proof of Theorem \ref{theorem 1}.\\
\end{proof}

\begin{proof}
	Proof of Corollary  \ref{theorem 2}. Let $\bar{g}=\psi^{-2}g$ be a conformal metric of $g$. We are assuming that $\psi(r)$
	and $h(r)$ are functions of $r$, where $r=\sum_{k=1}^{n}\varepsilon_{k}x_{k}^{2}$. Hence, we have
	
	\begin{align*}
	\psi_{,x_{i}}=2\varepsilon_{i}x_{i}\psi',\hspace{0,5cm} \psi_{,x_{i}x_{i}}=4x_{i}^{2}\psi''+2\varepsilon_{i}\psi', \hspace{0,5cm} \psi_{,x_{i}x_{j}}=4\varepsilon_{i}\varepsilon_{j}x_{i}x_{j}\psi''
	\end{align*}
	and 
	\begin{align*}
	h_{,x_{i}}=2\varepsilon_{i}x_{i}h',\hspace{0,5cm} h_{,x_{i}x_{i}}=4x_{i}^{2}h''+2\varepsilon_{i}h', \hspace{0,5cm} h_{,x_{i}x_{j}}=4\varepsilon_{i}\varepsilon_{j}x_{i}x_{j}h''.
	\end{align*}
	
	Substituting these expressions into \eqref{01}, we obtain
	\begin{align*}
	4\varepsilon_{i}\varepsilon_{j}(n-2)x_{i}x_{j}\psi''+4\varepsilon_{i}\varepsilon_{j}x_{i}x_{j}\psi h''+(2\varepsilon_{i}x_{i}\psi').(2\varepsilon_{j}x_{j}h')+(2\varepsilon_{j}x_{j}\psi').(2\varepsilon_{i}x_{i}f')=0,
	\end{align*}
	which is equivalent to 
	\begin{align*}
	4\varepsilon_{i}\varepsilon_{j}\left[(n-2)\psi''+\psi h''+2\psi' h'\right]x_{i}x_{j}=0.
	\end{align*}
	
	Since there exist $i\neq j$, such that $x_{i}x_{j}\neq 0$, we have
	\begin{align*}
	(n-2)\psi''+h''\psi+2\psi'h'=0.
	\end{align*}
	
	Similarly, considering the equation \eqref{2}, we obtain
	\begin{center}
		\begin{align*}
		4\psi\left[(n-2)\psi''+\psi h''+2\psi'h'\right]x_{i}^{2}+
		2\psi\left[\left(n-2\right)\psi'+\psi h'\right]\varepsilon_{i}+2\varepsilon_{i}n[1-2(n-1)\rho]\psi \psi'
		\end{align*}
		\begin{align*}
		+4\varepsilon_{i}\sum_{k=1}^{n}\varepsilon_{k}x_{k}^{2}\left\{\left(n-1\right)\left[\left(\rho n-1\right)(\psi')^{2}-2\rho\psi \psi''\right]-\psi \psi'h'+\psi \psi''\right\}=\lambda\varepsilon_{i}.
		\end{align*}
	\end{center}
	Note that $(n-2)\psi''+\psi h''+2\psi' h'=0$ and $r=\sum_{k=1}^{n}\varepsilon_{k}x_{k}^{2}$. Therefore, we obtain 
	\begin{center}
		\begin{align*}
		2\psi\left[\left(n-2\right)\psi'+\psi h'\right]+2n[1-2(n-1)\rho]\psi \psi'
		\end{align*}
		\begin{align*}
		+4r\left\{\left(n-1\right)\left[\left(\rho n-1\right)(\psi')^{2}-2\rho\psi \psi''\right]-\psi \psi'h'+\psi \psi''\right\}=\lambda.
		\end{align*}
	\end{center}
	This concludes the proof of Corollary \ref{theorem 2}.\\
\end{proof}

\begin{proof} Proof of the Theorem \ref{coro1}  It is well known (see, e.g., \cite{RO} or \cite{MP}) that if $\bar{g}=\frac{g}{\psi^{2}}$, then

	\begin{equation*}
		K_{\bar{g}}=(n-1)\left(2\psi\Delta_{g}\psi-n|\nabla_{g}\psi|^{2}\right).
	\end{equation*}
	
	How we are assuming that $\psi(r)$
	is a  functions of $r$, where $r=\sum_{k=1}^{n}\varepsilon_{k}x_{k}^{2}$, then we have that $K_{\bar{g}}=0$ if, and only, if

	\begin{equation*}
	-n\psi\psi'-2r\psi \psi''+nr\left(\psi'\right)^2=0, 
	\end{equation*}
which is equivalent to

	\begin{equation*}
	-\frac{n}{2r}\frac{\psi'}{\psi}+\frac{n}{2}\left(\frac{\psi'}{\psi}\right)^2-\frac{\psi''}{\psi}=0.
	\end{equation*}
	
	By equality $\frac{\psi''}{\psi}=\left(\frac{\psi'}{\psi}\right)'+\left(\frac{\psi'}{\psi}\right)^2$, follows that 
	\begin{equation*}
	-\frac{n}{2r}\frac{\psi'}{\psi}+\frac{n}{2}\left(\frac{\psi'}{\psi}\right)^2-\left(\frac{\psi'}{\psi}\right)'-\left(\frac{\psi'}{\psi}\right)^2=0.
	\end{equation*}
	
	Taking $y=\frac{\psi'}{\psi}$, the previous equation becomes
	
	\begin{equation}\label{18}
	y'=-\frac{n}{2r}y+\frac{n-2}{2}y^2.
	\end{equation}

	Note that equation \eqref{18} is an ordinary differential equation of Bernoulli. Therefore, you can determine all your solutions, whose general solution is given by
	\begin{equation}
	y^{-1}=Ce^{F}-\frac{(n-2)}{2}e^{F}\int e^{-F}dr,\hspace{0.5cm} \text{where}\hspace{0.5cm} F(r)=\frac{n}{2}\int
	\frac{1}{r}dr=\ln r^{\frac{n}{2}},
	\end{equation}
	C is an arbitrary constant (for more details see \cite{PO}). Thus 
	\begin{equation*}
	y^{-1}=C r^{\frac{n}{2}}- \frac{(n-2)}{2}r^{\frac{n}{2}}\int r^{-\frac{n}{2}}dr
	\end{equation*}
	equivalently, 
	\begin{equation*}
	y^{-1}=\left(C -\frac{(n-2)}{2}k_1\right)r^{\frac{n}{2}}+r,
	\end{equation*}
	where $k_1$ is a real number. Implies that 
	\begin{equation*}
	y^{-1}=Ar^{\frac{n}{2}}+r,
	\end{equation*}
	where $A=C -\frac{(n-2)}{2}k_1$. It follow that 
	\begin{equation*}
		y=\frac{r^{-\frac{n}{2}}}{A+r^{\frac{2-n}{2}}},
	\end{equation*}
	how $y=\frac{\psi'}{\psi}$, we get 
	\begin{equation}\label{19}
	\psi(r)=exp\left\{{\int\frac{r^{-\frac{n}{2}}}{A+r^{\frac{2-n}{2}}}dr}+\ln k_2\right\}
	\end{equation} 
	where $k_2\in \mathbb{R}_{+}^*$. Note that 
	\begin{equation}\label{20}
	\int\frac{r^{-\frac{n}{2}}}{A+r^{\frac{2-n}{2}}}dr=\ln\left(A+r^{\frac{2-n}{2}}\right)^{\frac{2}{2-2}},
	\end{equation}
	combining the equations \eqref{19} and \eqref{20}, the following that
	\begin{equation}\label{21}
	\psi(r)=k_2B^{\frac{2}{2-n}},
	\end{equation}
	where $B=A+r^{\frac{2-n}{2}}$. How $n \geq 3$ we obtain that 
	\begin{equation}\label{22}
	\psi(r)= \frac{k_2r}{\left(1+Ar^{\frac{n-2}{2}}\right)^{\frac{2}{n-2}}},
	\end{equation}
\end{proof}

\begin{proof} Proof of the Corollary \ref{coro2} If $K_{\bar{g}}=0$, by Theorem \ref{coro1} we get $\psi(r)= \frac{k_2r}{\left(1+Ar^{\frac{n-2}{2}}\right)^{\frac{2}{n-2}}}$. We will show that $\bar{g}=\frac{g_{0}}{\psi^2}$ is complete.

Consider the manifolds $M=\left(\mathbb{R}^n\setminus\{0\},\bar{g}=\frac{g_0}{\psi^2}\right)$, where $\psi(r)=\frac{k_2r}{\left(1+Ar^\frac{n-2}{2}\right)^\frac{2}{n-2}}$, $k_2\in\mathbb{R}_{+}^{*}$ and $N=\left(\mathbb{R}^n\setminus\{0\},g=\frac{g_0}{\varphi^2}\right)$, where  $\varphi(r)=\sqrt{r}$, and $g_0$ is a Euclidean metric. Note that 
\begin{equation*}
|v|_{\bar{g}}=\frac{1}{\psi}|v|_{g_0}\hspace{0.5cm} \textit{and} \hspace{0.5cm} |v|_{g}=\frac{1}{\varphi}|v|_{g_0}
\end{equation*}

By other hand, we get 
\begin{equation*}
|v|_{\bar{g}}=\frac{\left(1+Ar^\frac{n-2}{2}\right)^\frac{2}{n-2}}{k_2r}|v|_{g_0}=\frac{\left(1+Ar^\frac{n-2}{2}\right)^\frac{2}{n-2}}{k_2\sqrt{r}}\frac{1}{\sqrt{r}}|v|_{g_0},
\end{equation*}
thus, 
\begin{equation*}
|v|_{\bar{g}}=f(r)|v|_g,
\end{equation*}
where $f(r)=\frac{\left(1+Ar^\frac{n-2}{2}\right)^\frac{2}{n-2}}{k_2\sqrt{r}}$.

To find $c> 0$ such that $|v|_{\bar{g}}\geq c|v|_g$, just solve the following problem
\begin{equation*}
\min\limits_{r\in\mathbb{R_{+}^{*}}} f(r) 
\end{equation*}
The first derivative of $f$ takes us 
\begin{equation*}
f'(r)=\frac{r^{\frac{1}{2}}\frac{2}{n-2}\left(1+Ar^\frac{n-2}{2}\right)^\frac{4-n}{n-2}A\frac{n-2}{2}r^{\frac{n-4}{2}}-\frac{1}{2}r^{-\frac{1}{2}}\left(1+Ar^\frac{n-2}{2}\right)^\frac{2}{n-2}}{k_2r},
\end{equation*}
equaivalently,
\begin{equation*}
f'(r)=\frac{\left(1+Ar^\frac{n-2}{2}\right)^\frac{2}{n-2}}{k_2r}\left[\left(1+Ar^\frac{n-2}{2}\right)^{-1}r^{\frac{n-3}{2}}-\frac{1}{2r^{\frac{1}{2}}}\right].
\end{equation*}
Therefore,
\begin{equation*}
f'(r)=\frac{\left(1+Ar^\frac{n-2}{2}\right)^\frac{4-n}{n-2}}{k_2r^{\frac{3}{2}}}\left(Ar^{\frac{n-2}{2}}-1\right).
\end{equation*}

Given $f$ is a real function, we have that $r$ is a critical point if, and only if, $f'(r)=0$. How $r>0$, the minimum point candidate is given by 
\begin{equation*}
r=\frac{1}{A^{\frac{2}{n-2}}}.
\end{equation*}

Let's calculate the second derivative of $f$ and evaluate at this point, this is,
\begin{equation*}
f''(r)=\frac{1}{2k_2}\left[\frac{4-n}{n-2}\left(1+Ar^\frac{n-2}{2}\right)^\frac{6-2n}{n-2}A\frac{n-2}{2}r^{\frac{n-4}{2}}\left(Ar^{\frac{n-5}{2}}-r^{-\frac{3}{2}}\right)\right]
\end{equation*}
\begin{equation*}
+\frac{1}{2k_2}\left[\left(1+Ar^\frac{n-2}{2}\right)^\frac{4-n}{n-2}\left(\frac{n-5}{2}Ar^{\frac{n-7}{2}}+\frac{3}{2}r^{-\frac{5}{2}}\right)\right]
\end{equation*}
equivalently,
\begin{equation*}
f''(r)=\frac{1}{2k_2}\left(1+Ar^\frac{n-2}{2}\right)^\frac{6-2n}{n-2}\left(\frac{4-n}{2}A^2r^{\frac{2n-9}{2}}-\frac{(4-n)A}{2}r^{\frac{n-7}{2}}\right)
\end{equation*}
\begin{equation*}
+\frac{1}{2k_2}\left(1+Ar^\frac{n-2}{2}\right)^\frac{4-n}{n-2}\left(\frac{n-5}{2}Ar^{\frac{n-7}{2}}+\frac{3A}{2}r^{-\frac{5}{2}}\right)
\end{equation*}
implies that 
\begin{equation*}
f''(r)=-\frac{\left(1+Ar^\frac{n-2}{2}\right)^\frac{6-2n}{n-2}}{4k_2r^{\frac{5}{2}}}\left(
A^2r^{n-2}+6Ar^{\frac{n-2}{2}}-2Anr^{\frac{n-2}{2}}-3\right).
\end{equation*}

Now let's evaluate the second derivative at the point $r=\frac{1}{A^{\frac{2}{n-2}}}$, that is,
\begin{equation*}
f''\left(\frac{1}{A^{\frac{2}{n-2}}}\right)=-\frac{\left(1+A\left(\frac{1}{A^{\frac{2}{n-2}}}\right)^\frac{n-2}{2}\right)^\frac{6-2n}{n-2}}{4k_2\left(\frac{1}{A^{\frac{2}{n-2}}}\right)^{\frac{5}{2}}}\left(
A^2\left(\frac{1}{A^{\frac{2}{n-2}}}\right)^{n-2}+6A\left(\frac{1}{A^{\frac{2}{n-2}}}\right)^{\frac{n-2}{2}}-2An\left(\frac{1}{A^{\frac{2}{n-2}}}\right)^{\frac{n-2}{2}}-3\right),
\end{equation*} 
equivalently,
\begin{equation*}
f''\left(\frac{1}{A^{\frac{2}{n-2}}}\right)=-\frac{2^\frac{6-2n}{n-2}}{4k_2\frac{1}{A^{\frac{5}{n-2}}}}\left(1+6-2n-3\right),
\end{equation*} 
implies that 
\begin{equation*}
f''\left(\frac{1}{A^{\frac{2}{n-2}}}\right)=-\frac{2^\frac{10-4n}{n-2}}{k_2}A^{\frac{5}{n-2}}\left(4-2n\right).
\end{equation*} 

Therefore,
\begin{equation*}
f''\left(\frac{1}{A^{\frac{2}{n-2}}}\right)=\frac{2^\frac{8-3n}{n-2}}{k_2}A^{\frac{5}{n-2}}\left(n-2\right).
\end{equation*} 

How $n\geq3$, $A,k_2\in \mathbb{R}_{+}^{*}$, we get $f''\left(\frac{1}{A^{\frac{2}{n-2}}}\right)>0$, consequently $r=\frac{1}{A^{\frac{2}{n-2}}}$ it's a minimum point. Therefore,
\begin{equation}
f\left(\frac{1}{A^{\frac{2}{n-2}}}\right)=\frac{\left(1+A\left(\frac{1}{A^{\frac{2}{n-2}}}\right)^\frac{n-2}{2}\right)^\frac{2}{n-2}}{k_2\sqrt{\left(\frac{1}{A^{\frac{2}{n-2}}}\right)}}=\frac{\left(4A\right)^{\frac{1}{n-2}}}{k_2}
\end{equation}

Just take $c=\frac{\left(4A\right)^{\frac{1}{n-2}}}{k_2}$, we get $f(r)\geq c$, $\forall r$. Thus $|v|_{\bar{g}}\geq c|v|_{g}$, how $N=\left(\mathbb{R}^n\setminus\{0\},g=\frac{g_0}{\varphi^2}\right)$ is complete, implies that  $M=\left(\mathbb{R}^n\setminus\{0\},\bar{g}=\frac{g_0}{\psi^2}\right)$ is complete. Therefore, the proof is done.	
\end{proof}	

\begin{proof}Proof of the Theorem \ref{theorem 5} How $\left(\mathbb{R}^{n},\bar{g}\right)$ is a $\rho$--Einstein sotiton, with zero scalar curvature follows  by Theorem \ref{coro1} that 
	\begin{equation}\label{5}
	\psi(r)=k_2B^{\frac{2}{2-n}},
	\end{equation}
	where $B=A+r^{\frac{2-n}{2}}$.
	Consequently,
	\begin{equation}\label{6}
	\psi'(r)=k_2B^{\frac{2}{2-n}}r^{-\frac{n}{2}},\hspace{0.5cm}\psi''(r)=\frac{k_2n}{2}\left(B^{\frac{2}{2-n}}r^{-n}-B^{\frac{2}{2-n}}r^{-\frac{(n+2)}{2}}\right).
	\end{equation}

	Replacing the expressions found in \eqref{5} and \eqref{6} in \eqref{4},  we have
	\begin{equation*}
	2\left(n-1\right)\left(1-n\rho\right)k_2^2B^{\frac{2}{2-n}}B^{\frac{n}{2-n}}r^{-\frac{n}{2}}+2\left(n-1\right)\left(n\rho-1\right)k_2^2rB^{\frac{2n}{2-n}}r^{-n}
	\end{equation*}
	\begin{equation*}
	+2\left(1-2\left(n-1\right)\rho\right)rk_2B^{\frac{2}{2-n}}\frac{k_2n}{2}\left(B^{\frac{2(n-1)}{2-n}}r^{-n}-B^{\frac{n}{2-n}}r^{-\frac{(n+2)}{2}}\right)
	\end{equation*}
	\begin{equation*}
	+k_2^2B^{\frac{2}{2-n}}\left(B^{\frac{2}{2-n}}-2B^{\frac{n}{2-n}}r^{\frac{2-n}{2}}\right)h'=\frac{\lambda}{2},
	\end{equation*}	
	equivalently,
	\begin{equation*}
	2\left(n-1\right)\left(1-n\rho\right)k_2^2B^{\frac{2}{2-n}}B^{\frac{n}{2-n}}r^{-\frac{n}{2}}+2\left(1-2\left(n-1\right)\rho\right)\frac{k_2^2n}{2}B^{\frac{2}{2-n}}B^{\frac{2\left(n-1\right)}{2-n}}r^{1-n}
	\end{equation*}
	\begin{equation*}
	-2\left(1-2\left(n-1\right)\rho\right)\frac{k_2^2n}{2}B^{\frac{2}{2-n}}B^{\frac{n}{2-n}}r^{-\frac{n}{2}}+2\left(n-1\right)\left(n\rho-1\right)k_2^2B^{\frac{2\left(n-1\right)}{2-n}}B^{\frac{2}{2-n}}r^{1-n}
	\end{equation*}
	\begin{equation*}
	+k_2^2B^{\frac{2}{2-n}}\left(B^{\frac{2}{2-n}}-2B^{\frac{n}{2-n}}r^{\frac{2-n}{2}}\right)h'=\frac{\lambda}{2}.
	\end{equation*}
	How $B\neq0$, we get 
	
	\begin{equation*}
	2\left(n-1\right)\left(1-n\rho\right)B^{\frac{n}{2-n}}r^{-\frac{n}{2}}+\left(1-2\left(n-1\right)\rho\right)nB^{\frac{2\left(n-1\right)}{2-n}}r^{1-n}-\left(1-2\left(n-1\right)\rho\right)nB^{\frac{n}{2-n}}r^{-\frac{n}{2}}
	\end{equation*}
	\begin{equation*}
	-2\left(n-1\right)\left(1-n\rho\right)B^{\frac{2\left(n-1\right)}{2-n}}r^{1-n}+k_2^2B^{\frac{2}{2-n}}\left(B^{\frac{2}{2-n}}-2B^{\frac{n}{2-n}}r^{\frac{2-n}{2}}\right)h'=\frac{\lambda}{2k_2^2B^{\frac{2}{2-n}}}.
	\end{equation*}
	Consequently, 
	\begin{equation*}
	(n-2)\left(B^{\frac{n}{2-n}}r^{-\frac{n}{2}}-B^{\frac{2(n-1)}{2-n}}r^{1-n}\right)+\left(B^{\frac{2}{2-n}}-2B^{\frac{n}{2-n}}r^{\frac{2-n}{2}}\right)h'=\frac{\lambda}{2k_2^2B^{\frac{2}{2-n}}}
	\end{equation*}

	equivalently,
	\begin{equation*}
	2k_2^2B^{\frac{2}{2-n}}\left(n-2\right)\left(B^{\frac{n}{2-n}}r^{-\frac{n}{2}}-B^{\frac{2\left(n-1\right)}{2-n}}r^{1-n}\right)+2k_2^2B^{\frac{2}{2-n}}\left(B^{\frac{2}{2-n}}-2B^{\frac{n}{2-n}}r^{\frac{2-n}{2}}\right)h'=\lambda.
	\end{equation*}
	
	Note that $B^{\frac{2}{2-n}}-2B^{\frac{n}{2-n}}r^{\frac{2-n}{2}}\neq0$, otherwise $B=2r^{\frac{2-n}{2}}$ and 	Consequently $ B = 2A$.
	
	Thus,
	\begin{equation*}
	h'(r)=\frac{\lambda}{2k_2^2B^{\frac{2}{2-n}}\left(B^{\frac{2}{2-n}}-2B^{\frac{n}{2-n}}r^{\frac{2-n}{2}}\right)}+\left(2-n\right)\left(\frac{B^{\frac{n}{2-n}}r^{-\frac{n}{2}}-B^{\frac{2\left(n-1\right)}{2-n}}r^{1-n}}{B^{\frac{2}{2-n}}-2B^{\frac{n}{2-n}}r^{\frac{2-n}{2}}}\right).
	\end{equation*}
	
	Making 
	\begin{equation*}
	\varphi(r)=\frac{\lambda}{2k_2^2B^{\frac{2}{2-n}}\left(B^{\frac{2}{2-n}}-2B^{\frac{n}{2-n}}r^{\frac{2-n}{2}}\right)},\hspace{0.3cm}\text{and}\hspace{0.3cm} w(r)=\left(2-n\right)\left(\frac{B^{\frac{n}{2-n}}r^{-\frac{n}{2}}-B^{\frac{2\left(n-1\right)}{2-n}}r^{1-n}}{B^{\frac{2}{2-n}}-2B^{\frac{n}{2-n}}r^{\frac{2-n}{2}}}\right),
	\end{equation*}
	the first derivative of this equations leads us to 
	\begin{equation}
	\varphi'(r)=-\frac{\lambda}{2k_2^2}\frac{\left(nB^{-1}r^{-\frac{n}{2}}-(n+2)B^{-2}r^{1-n}\right)}{\left(B^{\frac{2}{2-n}}-2B^{\frac{n}{2-n}}r^{\frac{2-n}{2}}\right)^2},
	\end{equation}
	and
	\begin{equation}
	w'(r)=\frac{\left(2-n\right)}{2}\frac{\left(3nB^{\frac{2n}{2-n}}r^{-n}-4(n-1)B^{\frac{3n-2}{2-n}}r^{\frac{2-3n}{2}}+2(n-2)B^{\frac{4(n-1)}{2-n}}r^{2(1-n)}-nB^{\frac{n+2}{2-n}}r^{-\frac{(n+2)}{2}}\right)}{\left(B^{\frac{2}{2-n}}-2B^{\frac{n}{2-n}}r^{\frac{2-n}{2}}\right)^2}
	\end{equation}
	Replacing the functions in the  equation \eqref{3}, that is,
	\begin{equation*}
	(n-2)\psi''+\psi h''+2\psi'h'=0, 
	\end{equation*}
	we obtain 
	
	\begin{equation*}
	(n-2)\frac{k_2n}{2}\left(B^{\frac{2(n-1)}{2-n}}r^{-n}-B^{\frac{n}{2-n}}r^{-\frac{2+n}{2}}\right)\left(B^{\frac{4}{2-n}}-4B^{\frac{2+n}{2-n}}r^{\frac{2-n}{2}}+4B^{\frac{2n}{2-n}}r^{2-n}\right)
	\end{equation*}
	\begin{equation*}
	+\frac{(2-n)}{2}k_2B^{\frac{2}{2-n}}\left(3nB^{\frac{2n}{2-n}}r^{-n}-4\left(n-1\right)B^{\frac{2-3n}{2-n}}r^{\frac{2-3n}{2}}+2\left(n-2\right)B^{\frac{4\left(n-1\right)}{2-n}}r^{2\left(1-n\right)}-nB^{\frac{n+2}{2-n}}r^{-\frac{n+2}{2}}\right)
	\end{equation*}
	\begin{equation*}
	2\left(2-n\right)k_2B^{\frac{n}{2-n}}r^{-\frac{n}{2}}\left(B^{\frac{n}{2-n}}r^{-\frac{n}{2}}-B^{\frac{2\left(n-1\right)}{2-n}}r^{1-n}\right)\left(B^{\frac{2}{2-n}}-2B^{\frac{n}{2-n}}r^{\frac{2-n}{2}}\right)
	\end{equation*}
	\begin{equation*}
	-\frac{\lambda}{2k_2}\left(nB^{\frac{n}{2-n}}r^{-\frac{n}{2}}-(n+2)B^{\frac{2(n-1)}{2-n}}r^{1-n}\right)=0,
	\end{equation*}

	Which is equivalent to

	\begin{equation*}
	B^{\frac{2(n+1)}{2-n}}r^{-n}-2B^{\frac{3n}{2-n}}r^ {\frac{2-3n}{2}}+B^{\frac{2(2n-1)}{2-n}}r^{2(1-n)}=\frac{\lambda}{(n-2)^2k_2^2}\left(nB^{\frac{n}{2-n}}r^{-\frac{n}{2}}-(n+2)B^{\frac{2(n-1)}{2-n}}r^{1-n}\right).
	\end{equation*}
	
	Note that $nB^{\frac{n}{2-n}}r^{-\frac{n}{2}}-(n+2)B^{\frac{2(n-1)}{2-n}}r^{1-n}\neq0$, because otherwise $B=\frac{n+2}{n}r^{\frac{2-n}{2}}$. On the other hand,
	\begin{equation*}
	B^{\frac{2(n+1)}{2-n}}r^{-n}-2B^{\frac{3n}{2-n}}r^{\frac{2-3n}{2}}+B^{\frac{2(2n-1)}{2-n}}r^{2(1-n)}=0
	\end{equation*}
	therefore,
	\begin{equation*}
	\left(B^{\frac{n+1}{2-n}} r^{-\frac{n}{2}}-B^{\frac{2n-1}{2-n}} r^{1-n}\right)^2=0, 
	\end{equation*}
	and consequently $B=r^{\frac{2-n}{2}}$, but this is a contradiction, because $B=\frac{n+2}{n}r^{\frac{2-n}{2}}$. Therefore,
	
	\begin{equation}\label{lambda}
	\lambda=(n-2)^2k_2^2\frac{\left(B^{\frac{n+1}{2-n}} r^{-\frac{n}{2}}-B^{\frac{2n-1}{2-n}} r^{1-n}\right)^2}{\left(nB^{\frac{n}{2-n}} r^{-\frac{n}{2}}-(n+2)B^{\frac{2(n-1)}{2-n}} r^{1-n}\right)}.
	\end{equation}
	
	How $\lambda$ is constant, we have that 
	
	\begin{equation}
	\frac{d}{dr}\frac{\left(B^{\frac{n+1}{2-n}} r^{-\frac{n}{2}}-B^{\frac{2n-1}{2-n}} r^{1-n}\right)^2}{\left(nB^{\frac{n}{2-n}} r^{-\frac{n}{2}}-(n+2)B^{\frac{2(n-1)}{2-n}} r^{1-n}\right)}=0.
	\end{equation}
	
	if, and only if, 
	\begin{equation*}
	(4n-1)B^{\frac{4n}{2-n}} r^{-2n}-(4n-1)(n+2)B^{\frac{5n-2}{2-n}} r^{\frac{2-5n}{2}}-n^2B^{\frac{3n+2}{2-n}} r^{-\frac{3n+2}{2}}+(n+2)nB^{\frac{4n}{2-n}} r^{-2n}
	\end{equation*}
	\begin{equation*}
	+(2-5n)nB^{\frac{5n-2}{2-n}} r^{\frac{2-5n}{2}}-(2-5n)(n+2)B^{\frac{2(3n-2)}{2-n}} r^{2-3n}+(2n-1)nB^{\frac{2(3n-2)}{2-n}} r^{2-3n}
	\end{equation*}
	\begin{equation*}
	-\frac{1}{2}\left(3n^2+2n-4\right)B^{\frac{4n}{2-n}} r^{-2n}-(2n-1)(n+2)B^{\frac{7n-6}{2-n}} r^{\frac{6-7n}{2}}+(n^2+n-2)B^{\frac{5n-2}{2-n}} r^{\frac{2-5n}{2}}
	\end{equation*}
	\begin{equation*}
	+\frac{n^2}{2}B^{\frac{3n+2}{2-n}} r^{-\frac{3n+2}{2}}+(3n^2+2n-4)+\frac{n^2}{2}B^{\frac{5n-2}{2-n}} r^{-\frac{2-5n}{2}}-2(n^2+n-2)B^{\frac{2(3n-2)}{2-n}} r^{2-3n}
	\end{equation*}
	\begin{equation*}
	-n^2B^{\frac{4n}{2-n}} r^{-2n}-\frac{1}{2}\left(3n^2+2n-4\right)B^{\frac{2(n-2)}{2-n}} r^{2-3n}+\frac{n^2}{2}B^{\frac{5n-2}{2-n}} r^{-\frac{2-5n}{2}}+(n^2+n-2)B^{\frac{7n-6}{2-n}} r^{-\frac{6-7n}{2}}=0,
	\end{equation*}
	if, and only if,
	
	\begin{equation}\label{derivada}
	A^2\left[-n^2r^{\frac{n-2}{2}}A^2+(n^2+4)A+4(1-n)r^{\frac{2-n}{2}}\right]=0.
	\end{equation}
	
	We will prove that equation \eqref{derivada} is satisfied if, and only, if $A=0$. For this, consider 
	\begin{equation*}
	f(r)= -n^2r^{\frac{n-2}{2}}A^2+(n^2+4)A+4(1-n)r^{\frac{2-n}{2}}.
	\end{equation*}	
	If $f(r)= 0 \;\; \forall r>0$ and $A \ne 0$,   then its derivative also is zero. Since there is a single value of $r$ such that $f'(r) =0$ given by $ r= \left( \frac{4(n-1)}{n^2A^2}\right)^{\frac{1}{n-2}}$, we get a contradiction.
	
	Therefore the equation \eqref{derivada} is satisfied if, and only, if $A=0$. In this case $B=r^{\frac{2-n}{2}}$ .  Substituindo  $B=r^{\frac{2-n}{2}}$ in \eqref{lambda} we obtain that $ \lambda = 0$.  Therefore the proof is done.
\end{proof}

\begin{proof} Proof of the Theorem \ref{theorem 6}
	Follows by  Theorem \ref{theorem 5} that $\left(\mathbb{R}^{n},\bar{g}\right)$, $\bar{g}=\frac{1}{\psi^{2}}g$ is  a  $\rho$--Einstein sotiton  with scalar curvature $K_{\bar{g}}=0$
	is steady if, and only if, $A=0$. Besides that $ h(r)$ is constant. How $A=0$ we obtain  from Lemma \ref{coro1}, that $\psi(r)=k_2r$.
	Follow of the \cite{MA} that $\left(\mathbb{R}^{n},\bar{g}\right)$ have sectional curvature  zero. Therefore, we conclude que $\left(\mathbb{R}^{n},\bar{g}\right)$ is flat. The reciprocal is automatically satisfied.
\end{proof}

\begin{proof} Proof of the Corollary \ref{theorem 7}
	How $\left(\mathbb{M}^{n},\bar{g}\right)$ is locally conformally flat and rotationally symmetric, then locally the metric $\bar{g}$ is given by  $\bar{g}=\frac{1}{\psi^{2}}g$ where  $ \psi =\psi(r)$ and 
	$r=\sum\limits_{k=1}^nx_{k}^2$ and $ g$ is the euclidean metric.Therefore the results obtained in Theorem \ref{theorem 5} and Corollary \ref{theorem 6} are satisfied.
\end{proof}

\begin{proof} Proof of the Proposition \ref{prop} Consider the manifolds $M=\left(\mathbb{R}^n\setminus\{0\},g=\frac{g_0}{\varphi^2}\right)$, where $\varphi(r)=re^{-\left(1+r^{\frac{n-2}{2}}\right)^{\frac{2}{n-2}}}$, $k_2\in\mathbb{R}_{+}^{*}$ and $N=\left(\mathbb{R}^n\setminus\{0\},\bar{g}=\frac{g_0}{\psi^2}\right)$, where  $\psi(r)=\frac{k_2r}{\left(1+r^\frac{n-2}{2}\right)^\frac{2}{n-2}}$, and $g_0$ is a Euclidean metric. Note that 
	\begin{equation*}
	|v|_g=\frac{1}{\varphi}|v|_{g_0}\hspace{0.5cm} \textit{and} \hspace{0.5cm} |v|_{\bar{g}}=\frac{1}{\psi}|v|_{g_0}
	\end{equation*}
	
	By other hand, we get 
	\begin{equation*}
	|v|_{g}=\frac{e^{\left(1+r^{\frac{n-2}{2}}\right)^{\frac{2}{n-2}}}}{r}|v|_{g_0}=\frac{e^{\left(1+r^{\frac{n-2}{2}}\right)^{\frac{2}{n-2}}}}{\left(1+r^{\frac{n-2}{2}}\right)^{\frac{2}{n-2}}}\frac{\left(1+r^\frac{n-2}{2}\right)^\frac{2}{n-2}}{k_2r}|v|_{g_0}.
	\end{equation*}
	Thus, 
	\begin{equation*}
	|v|_{g}=h(r)|v|_{\bar{g}},
	\end{equation*}
	where, $h(r)=\frac{e^{\left(1+r^{\frac{n-2}{2}}\right)^{\frac{2}{n-2}}}}{\left(1+r^{\frac{n-2}{2}}\right)^{\frac{2}{n-2}}}$.
	
	To find $c_1> 0$ such that $|v|_{\bar{g}}\geq c_1|v|_g$, just solve the following problem
	\begin{equation*}
	\min\limits_{r\in\mathbb{R_{+}^{*}}} h(r) 
	\end{equation*}
	The first derivative of $h$ takes us 
	\begin{equation*}
	h'(r)=e^{\left(1+r^{\frac{n-2}{2}}\right)^{\frac{2}{n-2}}}\left(1+r^{\frac{n-2}{2}}\right)^{\frac{-n}{n-2}}r^{\frac{n-4}{2}}\left[\left(1+r^{\frac{n-2}{2}}\right)^{\frac{2}{n-2}}-1\right].
	\end{equation*}
	
	How $r>0$, we get $h'(r)>0$, $\forall r$, so the $h$ function is strictly increasing. Therefore, 
	\begin{equation}
	\min\limits_{r\in\mathbb{R}^{*}_+}h(r)=\lim\limits_{r\longrightarrow 0}\frac{e^{\left(1+r^{\frac{n-2}{2}}\right)^{\frac{2}{n-2}}}}{\left(1+r^{\frac{n-2}{2}}\right)^{\frac{2}{n-2}}}=e.
	\end{equation}
	
	Just take $c_1=e$, we get $h(r)\geq c_1$, $\forall r$. Thus $|v|_{g}\geq c|v|_{\bar{g}}$, how  $M=\left(\mathbb{R}^n\setminus\{0\}, \bar{g}=\frac{g_0}{\psi^2}\right)$ is complete, implies that $N=\left(\mathbb{R}^n\setminus\{0\}, g=\frac{g_0}{\varphi^2}\right)$ is complete. 
	
	We will show that $\left(\mathbb{R}^n,\bar{g}\right)$ has negative scalar curvature. Indeed, note that
	\begin{equation}\label{pride}
	\varphi'(r)=\frac{1-r^{\frac{n-2}{2}}\left(1+r^{\frac{n-2}{2}}\right)^{\frac{4-n}{n-2}}}{e^{\left(1+r^{\frac{n-2}{2}}\right)^{\frac{2}{n-2}}}},
	\end{equation}
and
\begin{equation*}
\varphi''(r)=-\frac{1}{2e^{\left(1+r^{\frac{n-2}{2}}\right)^{\frac{2}{n-2}}}}\left[(4-n)r^{\frac{2n-6}{2}}\left(1+r^{\frac{n-2}{2}}\right)^{\frac{6-2n}{n-2}}+(n-2)r^{\frac{n-4}{2}}\left(1+r^{\frac{n-2}{2}}\right)^{\frac{4-n}{n-2}}\right.
\end{equation*}
\begin{equation*}
\left.-2r^{\frac{n-4}{2}}\left(1+r^{\frac{n-2}{2}}\right)^{\frac{4-n}{n-2}}-2r^{\frac{2n-6}{2}}\left(1+r^{\frac{n-2}{2}}\right)^{2\frac{(4-n)}{n-2}}\right], 
\end{equation*}
implies that,
\begin{equation}\label{segder}
\varphi''(r)=-\frac{\left[nr^{\frac{n-4}{2}}\left(1+r^{\frac{n-2}{2}}\right)^{\frac{4-n}{n-2}}+(4-n)r^{n-3}\left(1+r^{\frac{n-2}{2}}\right)^{\frac{6-2n}{n-2}}-2r^{n-3}\left(1+r^{\frac{n-2}{2}}\right)^{2\frac{(4-n)}{n-2}}\right]}{2e^{\left(1+r^{\frac{n-2}{2}}\right)^{\frac{2}{n-2}}}}. 
\end{equation}

It is well known (see, e.g., \cite{RO} or \cite{MP}) that if $\bar{g}=\frac{g}{\psi^{2}}$, then 
\begin{equation}\label{curv}
K_{\bar{g}}=4r\left[2(n-1)\varphi\varphi''-n(n-1)\left(\varphi'\right)^2\right]+4n(n-1)\varphi\varphi'.
\end{equation}

Substituting the expressions found in \eqref{pride} and \eqref{segder} into \eqref{curv}, we 
\begin{equation*}
K_{\bar{g}}=\frac{4(n-1)}{e^{\left(1+r^{\frac{n-2}{2}}\right)^{\frac{2}{n-2}}}}\left[(2-n)r^{n-1}\left(1+r^{\frac{n-2}{2}}\right)^{\frac{2(4-n)}{n-2}}+2(1-n)r^{\frac{n}{2}}\left(1+r^{\frac{n-2}{2}}\right)^{\frac{4-n}{n-2}}\right.
\end{equation*}
\begin{equation*}
+\left.(n-4)r^{n-1}\left(1+r^{\frac{n-2}{2}}\right)^{\frac{2(3-n)}{n-2}}\right]
\end{equation*}
equivalently,
\begin{equation*}
K_{\bar{g}}=\frac{4(n-1)\left(1+r^{\frac{n-2}{2}}\right)^{\frac{4-n}{n-2}}}{e^{\left(1+r^{\frac{n-2}{2}}\right)^{\frac{2}{n-2}}}}\left[(2-n)r^{n-1}\left(1+r^{\frac{n-2}{2}}\right)^{\frac{4-n}{n-2}}+2(1-n)r^{\frac{n}{2}}\right.
\end{equation*}
\begin{equation*}
+\left.(n-4)r^{n-1}\left(1+r^{\frac{n-2}{2}}\right)^{-1}\right]
\end{equation*}
implies that,
\begin{equation*}
K_{\bar{g}}=\frac{4(n-1)\left(1+r^{\frac{n-2}{2}}\right)^{\frac{2(3-n)}{n-2}}}{e^{\left(1+r^{\frac{n-2}{2}}\right)^{\frac{2}{n-2}}}}\left[(2-n)r^{n-1}\left(1+r^{\frac{n-2}{2}}\right)^{\frac{2}{n-2}}+2(1-n)r^{\frac{n}{2}}\left(1+r^{\frac{n-2}{2}}\right)\right.
\end{equation*}
\begin{equation*}
-\left.(4-n)r^{n-1}\right].
\end{equation*}
Therefore,
\begin{equation*}
K_{\bar{g}}=-\frac{4(n-1)\left(1+r^{\frac{n-2}{2}}\right)^{\frac{2(3-n)}{n-2}}}{e^{\left(1+r^{\frac{n-2}{2}}\right)^{\frac{2}{n-2}}}}\left[(n-2)\left(1+r^{\frac{n-2}{2}}\right)^{\frac{2}{n-2}}+2(n-1)r^{\frac{2-n}{2}}+(n+2)\right]r^{n-1}.
\end{equation*}

\end{proof} 
\begin{proof} Proof of Corollary \ref{coro3} 
Follow immediately from Corollary \ref{coro2}.
\end{proof}

The authors would like to thank the referee for his careful reading, relevant remarks, and valuable suggestion.

\end{document}